\documentclass[12pt,leqno,fleqn,epsfig]{article}
\usepackage{amssymb, epsfig, amsmath, amsthm}
\usepackage{mathrsfs}       

\newcommand{\R}{\mathbb{R}}

\newcommand{\bu}{\boldsymbol u}

\newcommand{\bL}{\boldsymbol L}

\newcommand{\bfvarphi}{\boldsymbol\varphi}

\newcommand{\const}{\operatorname*{const}}

\newcommand{\be}{\begin{equation}}
\newcommand{\ee}{\end{equation}}
\newcommand{\bea}{\begin{eqnarray}}
\newcommand{\eea}{\end{eqnarray}}
\newcommand{\bean}{\begin{eqnarray*}}
\newcommand{\eean}{\end{eqnarray*}}

\newcommand{\var}{\varepsilon}
\renewcommand{\o}{\omega}

\newcommand{\intl}{\int\limits}

\newcommand{\Beweisende}{\rule{0.2cm}{0.2cm}}

\newcounter{secnum}

\newtheorem{thm}{Theorem}[section]

\newtheorem{lem}[thm]{Lemma}

\theoremstyle{definition}

\newtheorem{defin}[thm]{Definition}
\newtheorem{rem}[thm]{Remark}

\title{On the Liouville theorem for  weak Beltrami flows } 
 
\author{Dongho Chae$^*$  and J\"{o}rg Wolf $^\dagger$\\
\ \\
 $*$Department of Mathematics\\
Chung-Ang University\\
 Seoul 156-756, Republic of Korea\\
 e-mail: dchae@cau.ac.kr\\
and \\
$\dagger$Department of Mathematics\\
Humboldt University Berlin\\
Unter den Linden 6, 10099 Berlin, Germany\\
e-mail: jwolf@math.hu-berlin.de}
\date{}
\begin{document}

\maketitle
\begin{abstract}
We study Beltrami flows in the setting of weak solution to the stationary Euler equations in $\Bbb R^3$. For this weak Beltrami flow we prove the regularity and the Liouville  property. In particular, we show that if  tangential part of the velocity has certain decay property at infinity, then the solution becomes trivial. This decay condition of of the velocity is weaker than the previously known sufficient conditions for the Liouville property of the Betrami flows.  For the proof we establish a mean value formula  and other various formula for the tangential and the normal components of the weak solutions to the stationary  Euler equations.\\
\ \\
\noindent{\bf AMS Subject Classification Number:}  76B03, 35Q31\\
  \noindent{\bf
keywords:} Euler's equations, Beltrami flows, Liouville's theorem

\end{abstract}

\section{Introduction}
\label{sec:-1}
\setcounter{secnum}{\value{section} \setcounter{equation}{0}
\renewcommand{\theequation}{\mbox{\arabic{secnum}.\arabic{equation}}}}

We consider the stationary Euler equation in $ \R^{3}$
\begin{equation}
(\bu \cdot \nabla ) \bu  = - \nabla p,\qquad  \nabla \cdot  \bu =0,
\label{1.1}
\end{equation}
where $\bu=(u_1 (x), u_2 (x), u_3 (x))$,  $x\in  \Bbb R^3$.  Together with the Navier-Stokes equations the Euler system
is of  fundamental importance  in the fluid mechanics. Even with such case the stationary Euler equations is less well studied than the Navier-Stokes equations.  On the other hand the Euler equations have  rich  geometric structure as  can be seen in the book by Arnold and Khesin(\cite{arn}) for example. 
Using the vector identity 
\be 
\label{1.1a}
\frac12 \nabla |\bu|^2 =(\bu \cdot \nabla )\bu +\bu \times (\nabla \times \bu),
\ee
 one can rewrite the system (\ref{1.1}) as
\be\label{1.1b}
\bu \times \omega=\nabla (p+\frac12 |\bu|^2), \quad \omega=\nabla \times \bu.
\ee
A vector field in $\Bbb R^3$, $\bu$ is called the Beltrami field if there exists a scalar function $\lambda=\lambda (x)$ such that $\o=\lambda \bu$.  If  $\bu$  is a Betrami field,  and simultaneously is a solution to (\ref{1.1}), then  the pressure should be given by $p=-\frac12 |\bu|^2$ up to addition of a constant. Conversely,  if $(\bu, p)$ solves (\ref{1.1}) with $p=-\frac12 |\bu|^2$, then $\bu$ is a Beltrami field. A Beltrami field, which is a solution to the stationary Euler system, is called the Beltrami flow.
The Beltrami fields has special roles in the physics of turbulence. In fact it is  observed numerically in \cite{pel} that 
 turbulent flow  consists of the superposition of the Beltrami fields. This aspect was theoretically studied in \cite{con}, and 
 more recently the Beltrami fields are crucially used to construct energy dissipating continuous weak solutions of the time dependent Euler equations(\cite{del}). For the study of Beltrami flows the authors of \cite{enc} constructed a non-trivial 
 solution $\bu(x)$  which decays of $O(\frac{1}{|x|})$ as $|x|\to +\infty$. Remarkably enough, in \cite{nad} it is shown later that if either the Beltrami flow $\bu$  decays  like $o(\frac{1}{|x|})$ as $|x|\to +\infty$ or $\bu \in L^q(\Bbb R^3)$, $q\in [2, 3]$, then $\bu\equiv 0$. Another sufficient condition
 for this Liouville type property $\frac{|\bu|^2}{|x|^{\mu}} \in L^1(\Bbb R^3)$, $\mu \in [0, 1]$ is also found in \cite{cha}.  These  studies on the Beltrami flow was in the setting of classical formulation, assuming $C^1(\Bbb R^3)$ regularity or higher for the velocity fields. 
  In this paper we formulate the notion of the Beltrami flow in the weak sense, and study its regularity as well as the Liouville type property. Based on the mean value formula as derived in the next section we could obtain, among others, much weaker condition for the Liouville property, imposing the decay condition only on the {\em tangential component} of the velocity.\\
 
We start from the weak formulation of (\ref{1.1}).
\begin{defin}
A pair $ (\bu , p)\in  [L^2_{ \rm loc}(\R^{3})]^3\times L^1_{ \rm loc}(\R^{3}) $ is said to be a {\it weak solution} to \eqref{1.1} 
if 
\begin{equation}
\intl_{ \R^{3}}  \bu \otimes \bu:  \nabla  \bfvarphi   dx =   - \intl_{ \R^{3}} p \nabla \cdot \bfvarphi   dx  \quad  
\forall\, \bfvarphi \in  [C^{\infty}_{\rm c}(\R^{3})]^3.
\label{1.1a}
\end{equation}
\end{defin}
Then, we can  generalize the classical notion of the Beltrami flow as follows.
\begin{defin} A  weak solution $ (\bu, p) $ to (\ref{1.1})   is called  weak  {\it Beltrami flow} if it satisfies
\begin{equation}
p = -\frac{| \bu |^2}{2} \quad  \text{ a.\,e. in}\quad  \R^{3}.
\label{1.1b}
\end{equation}
\end{defin}

\hspace{0.5cm}
We denote the normal and the tangential part of the velocity $\bu(x)$  are defined respectively as follows.
\[
\bu _N := \Big(\bu \cdot  \frac{x}{|x|}\Big) \frac{x}{|x|}, \qquad 
\bu _T := \bu \times  \frac{x}{|x|}, 
\]
where the equality holds almost everywhere in $\Bbb R^3$. Then,  we have the orthogonal decomposition
for any vector field on $\Bbb R^3$.
$$ \bu= \bu_N+\bu_T. $$
The main purpose of this paper  is to present  decay conditions  on $ \bu _T = \bu \times \frac{x}{|x|}$, 
which   imply the Liouville property $ \bu \equiv {\bf 0}$ for  the weak Beltrami flows.  This would extend the previous Liouville type  results in  \cite{nad, cha}.   In addition, we show that for a weak Beltrami flow the velocity in fact enjoys certain regularity properties in terms of weighted Lebesgue spaces.  Both results will be achieved from various identities involving the tangential part and the  normal part  of $ \bu $.  

We first state the following Liouville  property for a weak Beltrami flow.
\begin{thm}
\label{thm1.3}
Let $ \bu \in \bL^2_{ \rm loc}$ be a weak Beltrami flow.  Then 
\begin{equation}
 \intl_{B_R} \frac{| \bu_N |^2}{|x|}  dx \le  \frac{1}{2 R} \intl_{B_R}  | \bu |^2 
 dx  = \frac{1}{2} \intl_{\partial B_R} (| \bu_T|^2 - | \bu_N|^2) dS  \quad  \forall\, 0< R < +\infty.  
\label{2.12}
\end{equation}
Therefore, if there exists a sequence $ R_k \rightarrow  +\infty$ such that 
\begin{equation}
\intl_{\partial B_{ R_k}} | \bu_T |^2dS \rightarrow  0 \quad  \text{as}\quad  k \rightarrow +\infty
\label{2.13}
\end{equation}
then $ \bu \equiv 0$. 
 \end{thm}

 \begin{rem}
 As  immediate consequences of the above theorem we obtain the previously known Liouville type results in \cite{nad, cha}.
 Indeed, let $ \bu \in \bL^2_{ \rm loc}$ be a weak Beltrami flow, satisfying one of the following conditions.
 \begin{itemize}
 \item[(i)] $ |\bu_T (x)| =o(|x|^{-1})$ as $|x|\to +\infty$. 
\item[(ii)] $ \bu_T\in [L^q(\Bbb R^3)]^3$ for some $ q\in [2, 3]$. 
\item[(iii)] $ \frac{ |\bu_T |^2}{|x|^\mu }\in L^1(\R^{3})$ for some $\mu \in (-\infty, 1]$.
\end{itemize}
Then, one can find a sequence $\{ R_k\}_{k\in \Bbb N}$ such that (\ref{2.13}) holds. The case  (i)  is obvious. In the case (ii) from the condition
$$ \|\bu_T\|_{L^q} ^q=\int_0 ^{+\infty} \int_{\partial B_R} |\bu_T|^q dS dR<+\infty 
$$
 one can deduce  that  there exists a sequence $\{ R_k\}_{k\in \Bbb N}$ such that
$$ \int_{\partial B_{R_k}} |\bu_T|^q dS=o\left(\frac{1}{R_k}\right) \quad \mbox{as}\quad k\to +\infty.
$$
Therefore, for each $q\in (2, 3]$, we obtain 
$$
\int_{\partial B_{R_k}} |\bu_T|^2 dS \leq C\left( \int_{\partial B_{R_k}} |\bu_T|^q dS\right)^{\frac{2}{q}} R_k ^{\frac{2(q-2)}{q}} =
 o (R_k ^{-\frac{2}{q}} )  R_k ^{\frac{2(q-2)}{q}}=  o(1)\quad \mbox{as}\quad k\to +\infty.
$$
 In the case (iii) we have for $\mu\in (-\infty, 1]$
$$\int_0 ^{+\infty} R^{-\mu}\int_{\partial B_R}| \bu_T|^2dS dR  <+\infty, $$
from which we find there exists a sequence $\{ R_k\}_{k\in \Bbb N}$ such that
$$ \int_{\partial B_{R_k}} |\bu_T|^2 dS=o\left(\frac{1}{R_k}\right) R_k ^\mu=o(1)\quad \mbox{as}\quad k\to +\infty.
$$
\end{rem}

\begin{rem}
If $ \bu $ is a non-trivial Beltrami flow such that $ |\bu | \le  \frac{K}{|x|}$ $ (K=\const >0)$, as has been constructed 
in \cite{enc} for example, then 
\begin{equation}
\intl_{ \R^{3}} \frac{ | \bu_T |^2}{|x|}  dx = +\infty,\quad \text{and}\quad     
\intl_{ \R^{3}} \frac{ | \bu_N |^2}{|x|}  dx <+\infty.
\label{2.15}
\end{equation}
Hence, $ \bu _N$ has stronger decay then $ \bu _T$. 
Clearly the first statement follows directly from (iii) of Remark 1.4, while the second one follows from \eqref{2.12} together with 
$$\intl_{B_R} \frac{| \bu_N |^2}{|x|}  dx\le \frac{1}{2} \intl_{\partial B_R} (| \bu_T|^2 - | \bu_N|^2) dS \leq2\pi K^2 \quad  \forall\, 0< R < +\infty,$$
and passing $R\to +\infty.$
\end{rem} 

\vspace{0.5cm}  

Let us recall that we say a vector field $\bu\in [L^p_{\rm loc} (\Bbb R^3)]^3$, $1\leq p\leq +\infty$, belongs to the local Morrey space $[L^{p,\lambda}_{\rm loc} (\Bbb R^3)]^3$ if 
$$\sup_{\rho>0}\left( \rho ^{-\lambda} \int _{B_\rho } |\bu|^p dx\right) <+\infty. $$
The following is a regularity result for weak Beltrami flows. 
\begin{thm}
\label{thm1.6}
Let $ \bu \in \bL^2_{ \rm loc}$ be a weak Beltrami flow. Then,  $ \bu \in [L^{ 2,1}_{ \rm loc}(\R^{3})]^3$. 
\end{thm} 

\section{Mean value formulas for weak solutions to the stationary Euler equation}
\label{sec:-2}
\setcounter{secnum}{\value{section} \setcounter{equation}{0}
\renewcommand{\theequation}{\mbox{\arabic{secnum}.\arabic{equation}}}}

The aim of this section is to derive identities obtained from \eqref{1.1a}. We have the following 

\begin{lem}
\label{lem2.1}
Let $ (\bu ,p)$ be a weak solution to the stationary Euler equation. Then there holds 

\begin{align}
\intl_{\partial B_R} (p + | \bu_N |^2 ) dS
= \frac{1}{R}\intl_{B_R} (3p + | \bu |^2 ) dx  \quad  \text{for a.\,e. $ 0<R<+\infty$}.
\label{2.1}
\end{align}

\end{lem}
\begin{rem}
A formula similar to (\ref{2.1}) is derived in Theorem 1.2 of \cite{cha0}, where
the integration over $B_R$ is replaced by the integration over $\Bbb R^3\setminus B_R$. In the above lemma our assumptions on $(\bu, p)$  are much weaker in terms of the regularity and the integrability, and  the method of proof given below is completely different from \cite{cha0}.
\end{rem}
{\bf Proof of Lemma \ref{lem2.1} } We first regularize the Euler equations, using the standard mollifying  kernel 
$ \rho _\var = \var ^{ -3} \rho (\var ^{ -1}x)$, $ 1< \var <+\infty$,  where 
$\rho \in C^{\infty}_{\rm c}(B_1)$ is nonnegative, and radial symmetric together    
with the normalized condition 
\begin{equation}
\label{2.1a}
\intl_{ \R^{3}} \rho   dx=1.
\end{equation}
For $ f \in L^1_{ \rm loc}(\R^{3})$ we define
\[
f_\var (x) = \intl_{ \R^{3}} f(y) \rho_\var  (x-y)  dy,\quad  x \in \R^{3}.    
\]     
We note that \eqref{2.1a} implies
\begin{equation}
\intl_{ B_R} | f_\var  |  dx \le   \intl_{ B_{R+\var }} | f|  dx\quad  \forall\,0<R<+\infty.
\label{2.3}
\end{equation}
  
From \eqref{1.1a} with $ \bfvarphi_\var   $ in place of $ \bfvarphi $ we easily see that 
\begin{equation}
\nabla \cdot (\bu \otimes \bu)_\var  = - \nabla p_\var  \quad  \text{in}\quad  \R^{3}. 
\label{2.4}
\end{equation} 
We multiply both sides of \eqref{2.4} by $ x $, and for $ 0< r < +\infty$  we integrate the result over $ \partial B_r $ with respect to the two dimensional surface measure. This yields
\begin{align}
\intl_{\partial B_r }   \partial _j (u^j u^i)_\var   x_i   d S =- \intl_{\partial B_r} \partial _i p_\var   x_i  dS,
\label{2.4a}
\end{align}
where and hereafter we use the Einstein convention that the repeated indices imply the summation over $\{ 1,2,3\}$.
Let $-\infty< \alpha < 3$. We multiply both sides of this by $ r ^{ -\alpha } = |x|^{-\alpha }$, and  integrate 
it over $ (0,R)$, $ 0<R<+\infty$ with respect  to $ r$, and then apply the integration by parts. Thus, we have
\begin{align*}
\intl_{B_R}  \partial _j  (u^j u^i )_\var  x_i |x|^{ -\alpha }  d x 
&=- \intl_{B_R}    (|\bu|^2)_\var   |x|^{ -\alpha }  d x  + \alpha  \intl_{B_R}   (u^i u^j )_\var   x_ix_j  |x|^{ -\alpha -2}  d x
\\
&\qquad \qquad   + \intl_{\partial B_R}   (u^i u^j )_\var   x_ix_j|x|^{ -\alpha -1} dS  
\\
&=- \intl_{B_R} \partial _i p_\var   x_i |x|^{ -\alpha } dx  
\\
&= (3 -\alpha )\intl_{B_R} p_\var   |x|^{ -\alpha } dx - \intl_{\partial B_R} p_\var  |x|^{ -\alpha +1}  dS.  
\end{align*}
From the above identity after  a routine manipulation we obtain 
\begin{align}
R^{ 1-\alpha }
\intl_{\partial B_R} (p_\var   + v _{ N, \var }   ) dS=\intl_{B_R} \left\{(3-\alpha)p_\var   + (| \bu |^2)_\var   -
\alpha  v_{ N, \var  }\right\} |x|^{ -\alpha }  dx,  
\label{2.5}
\end{align}  
where we have set 
\begin{align}\label{2.5a}
v_{ N, \var  }(x) &= (u^i u^j )_\var   \frac{x_ix_j}{|x|^2} =  \frac{1}{|x|^2} \intl_{ \R^{3}} 
(\bu(y) \cdot x )^2 \rho_\var  (x-y)  dy\nonumber\\
&=\int_{\Bbb R^3} |\bu_N (x)|^2 \rho_\var (x-y)dy,  \quad  \forall x \in \R^{3}  \setminus  \{0\}.  
\end{align}
In particular,  \eqref{2.5} with $ \alpha =0$  yields
\begin{align}
\intl_{\partial B_r}( p_\var   + v _{ N, \var }   ) dS= \frac{1}{r}\intl_{B_r} 3p_\var   + (| \bu |^2)_\var 
 dx  \quad  \forall\, 0 \le r < +\infty.
\label{2.6}
\end{align}  
For an arbitrary  $ h>0$ we integrate the above identity over $ (R, R+h)$, and then letting $ \var \rightarrow 0$. This leads to 
\begin{equation}
\frac{1}{h} \intl_{R}^{R+h}\intl_{\partial B_r} ( p   + | \bu _N|^2)  dSdr=
\frac{1}{h}\intl_{B_{ R+h}\setminus  B_R} ( p   + | \bu _N|^2) dx= 
\frac{1}{h} \intl_{R}^{R+h}   
\frac{1}{r} \intl_{B_r}( 3p  + | \bu |^2 )
 dx dr.
\label{2.7}
\end{equation}

Since $ r \mapsto \intl_{\partial B_r} ( p   + | \bu _N|^2) dS$ is $ L^1_{ \rm loc}([0, +\infty))$, we have 
\[
\lim_{ h\to 0} \frac{1}{h} \intl_{R}^{R+h}\intl_{\partial B_r} ( p   + | \bu _N|^2 ) dSdr = \intl_{\partial B_R} ( p   + | \bu _N|^2 ) dS
\quad  \text{for a.\,e. $ 0<R<+\infty$}. 
\]
Hence, the assertion follows  from \eqref{2.7}  after letting $ h \rightarrow 0$ on both sides.  \hfill \Beweisende 

\begin{rem}
Define, 
\begin{equation}\label{2.7a}
\varphi (r) = -2\intl_{\partial B_r} ( p   + | \bu _N|^2 ) dS,\quad  0< r< +\infty.
\end{equation}
By virtue of  \eqref{2.7} we see that each $ r\in (0, +\infty)$ is a Lebesgue point of $ \varphi $. Therefore 
$  \varphi $ can be identified with the continuous function    $ r \mapsto -\frac{2}{r} \intl_{B_r} ( 3p  + | \bu |^2 )
 dx$, which is  actually  absolutely continuous.  Indeed, let us set $-2(3p+|\bu|^2 ):=g\in L^1_{\rm{loc}} (\Bbb R^3)$, and
 $\{ (a_i, b_i)\}_{i\in J}$ be a collection of non-overlapping intervals  such that $0<r=\inf_{i\in J}a_i,$ and $R=\sup_{i\in J} b_i<+\infty$. Then,
 \begin{align*}
\sum_{i\in J} | \varphi (b_i) -\varphi (a_i)|=\sum_{i\in J}\left| \frac{1}{b_i}\int_{B_{b_i}\setminus B_{a_i}} gdx -\frac{b_i-a_i}{a_i b_i }\int_{B_{a_i}} gdx \right| \\
\leq  \frac{1}{r} \sum_{i\in J} \int_{B_{b_i}\setminus B_{a_i}} |g|dx + \frac{1}{r^2} \sum_{i\in J} (b_i-a_i)\int_{B_{R}} |g|dx \to 0
 \end{align*}
 as $\sum_{i\in J} (b_i-a_i) \to 0$.
   By  the similar argument we have used in the proof of Lemma\,\ref{lem2.1},  from \eqref{2.5} we claim that  the following general identity   holds true.
\begin{equation}
R^{ 1-\alpha }\varphi (R) - r^{ 1-\alpha }\varphi (r) = 
-2\intl_{B_R  \setminus B_r } \left\{(3-\alpha)p   + | \bu |^2   -
\alpha  | \bu_{ N }|^2\right\} |x|^{ -\alpha }  dx\
\label{2.8}
\end{equation}
for all  $0<r <R$ and   $\alpha \in  \R$.  
Indeed, given $ 0< \sigma <\rho $ we multiply both sides of \eqref{2.4a} (with $ 0< \var < \sigma  $) by  $ r ^{ -\alpha } = |x|^{-\alpha }$. Then we integrate the result  over $ (\sigma, \rho )$, and  integrate by parts. This yields 
\begin{align*}
&\intl_{B_\rho \setminus B_\sigma  }  \partial _j  (u^j u^i )_\var  x_i |x|^{ -\alpha }  d x 
\\
&\quad =- \intl_{B_{ \rho }  \setminus B_\sigma }    (|\bu|^2)_\var   |x|^{ -\alpha }  d x  + \alpha  \intl_{B_\rho   \setminus 
B_\sigma }  
 (u^i u^j )_\var   x_i x_j  |x|^{ -\alpha -2}  d x
\\
&\qquad \qquad   + \intl_{\partial B_\rho }   (u^i u^j )_\var   x_ix_j|x|^{ -\alpha -1} dS- \intl_{\partial B_\sigma  }   (u^i u^j )_\var   x_ix_j|x|^{ -\alpha -1} dS  
\\
&\quad =- \intl_{B_\rho   \setminus B_\sigma} \partial _i p_\var   x_i |x|^{ -\alpha } dx  
\\
&\quad = (3 -\alpha )\intl_{B_\rho   \setminus B_\sigma  } p_\var   |x|^{ -\alpha } dx - \intl_{\partial B_\rho } p_\var  |x|^{ -\alpha +1}  dS
+ \intl_{\partial B_\sigma } p_\var  |x|^{ -\alpha +1}  dS.  
\end{align*}
 Accordingly, 
 \begin{align}
& \rho ^{ 1-\alpha }
\intl_{\partial B_\rho } (p_\var   + v _{ N, \var }   ) dS - \sigma  ^{ 1-\alpha }
\intl_{\partial B_\sigma } (p_\var   + v _{ N, \var }   ) dS 
 \cr
&\qquad   \qquad =\intl_{B_\rho \setminus B_\sigma } \left\{(3-\alpha)p_\var   + (| \bu |^2)_\var   -
\alpha  v_{ N, \var  }\right\} |x|^{ -\alpha }  dx. 
 \label{2.6a}
 \end{align}

Let $ 0<r < R$ be arbitrarily chosen and let $ 0<\var < r$. 
Given $0<h < r  -\var $, we intgrate  (\ref{2.6a}) first  over $(R, R+h)$ with respect to $ \rho $ and second over 
$ (r, r+h)$ with respect to $ \sigma $. Then dividing both sides by $ h^2$,  we  have
\begin{align}
\frac{1}{h} \int_R ^{R+h} & \rho ^{ 1-\alpha }
\int_{\partial B_\rho } (p_\var   + v _{ N, \var }   ) dSd\rho  - \frac{1}{h} \int_r ^{r+h}\sigma ^{ 1-\alpha }
\intl_{\partial B_\sigma  } (p_\var   + v _{ N, \var }   ) dS d\sigma 
 \nonumber\\
\label{2.8a}
  &=\frac{1}{h^2} \int_r ^{r+h}\int_R ^{R+h} \intl_{B_\rho   \setminus B_\sigma } \left\{(3-\alpha)p_\var   + (| \bu |^2)_\var   -
\alpha  v_{ N, \var  }\right\} |x|^{ -\alpha }  dx d\rho d \sigma,
\end{align}  
where we used the notation of (\ref{2.5a}). Passing $\var \to 0$, using the convergence properties of the mollified functions,
we obtain that the left hand side of (\ref{2.8a}) converges to 
\begin{align}
&\frac{1}{h} \int_R ^{R+h} \rho ^{ 1-\alpha }
\int_{\partial B_\rho  } (p  + |\bu_N|^2  ) dSd\rho  - \frac{1}{h} \int_r ^{r+h}\sigma  ^{ 1-\alpha }
\intl_{\partial B_\sigma  } (p   + | \bu|_{ N}   ) dS d\sigma 
\nonumber \\
&\qquad =  -\frac{1}{2h}\intl_{R}^{R+h}\rho ^{ 1-\alpha } \varphi (\rho ) d\rho 
+\frac{1}{2h}\intl_{r}^{r+h}\sigma ^{ 1-\alpha } \varphi (\sigma  ) d\sigma.   
\label{2.8b}
\end{align}
Passing $\var\to 0$ in the right-hand side of (\ref{2.8a}), we first observe that
$$
\intl_{B_\rho   \setminus B_\sigma } \left\{(3-\alpha)p_\var   + (| \bu |^2)_\var   -
\alpha  v_{ N, \var  }\right\} |x|^{ -\alpha }  dx\to \intl_{B_\rho    \setminus B_\sigma } \left\{(3-\alpha)p   + | \bu |^2   -
\alpha  |\bu_N|^2\right\} |x|^{ -\alpha }  dx
$$
 for each $\rho \in (R, R+h)$ and $ \sigma \in (r, r+h)$. Thus, by means of Lebesgue's theorem of 
dominated convergence we see that the right-hand side of (\ref{2.8a})  tends to  
\begin{equation}
\label{2.8c}
\frac{1}{h^2} \int_r ^{r+h}\int_R ^{R+h} \intl_{B_\rho    \setminus B_\sigma } \left\{(3-\alpha)p   + | \bu |^2   -
\alpha  |\bu_N|^2\right\} |x|^{ -\alpha }  dx d\rho d \sigma
\end{equation}
as  $ \var \rightarrow 0$. With the help the above limits we are able to carry out the   passage to the limit $ \var \rightarrow 0$ in  \eqref{2.8a} .  Letting 
$ h \rightarrow 0$ in both of (\ref{2.8b}) and (\ref{2.8c}) and equating them to each other, we obtain \eqref{2.8}. The claim is proved.\\

In the case when the mean value $ (3p + | \bu |^2)_{ B_r}$ is bounded for $ r \rightarrow 0^+$,  
 it follows that $ r^{ -2}\varphi (r)$ is bounded too as $  r \rightarrow  0^+$, which shows that 
 $ r^{ -1-\alpha } \varphi (r) 
 \rightarrow 0$ as $ r \rightarrow 0^+$ for all $ \alpha < 1$. However, this property might not be true in general. 
 
\hspace{0.5cm}
 Furthermore, from \eqref{2.8} we easily get $ \varphi \in W^{1,\, 1}_{ \rm loc}((0, +\infty))$, and there holds 
 \begin{equation}
 \varphi + R \varphi ' = -2 \intl_{\partial B_R}( 3p + | \bu |^2 ) dS\quad  \text{for a.\,e. $ 0<R<+\infty$},
 \label{2.9}
 \end{equation} 
 and, substituting $\varphi $ from (\ref{2.7a}) into (\ref{2.9}), using the fact $|\bu|^2 =|\bu_N|^2 +|\bu_T|^2$, we arrive at 
 \begin{equation}
 \varphi ' = - \frac{2}{R} \intl_{\partial B_R} (2p + | \bu_T |^2  ) dS
 = \frac{2}{R}\varphi + \frac{2}{R}\intl_{\partial B_R}( 3| \bu_N |^2 - | \bu |^2 ) dS
 \label{2.9a}
 \end{equation}
for a.\,e. $ 0<R<+\infty$. 
\end{rem}

\hspace{0.5cm}
 In case of Beltrami flow, we have  $ 2p =-| \bu |^2$, and we get from \eqref{2.1} 
 \begin{equation}
 \varphi (r) = \intl_{\partial B_r}( | \bu_T |^2  - | \bu _N|^2) dS = \frac{1}{r} \intl_{B_r}  | \bu |^2 
 dx\quad  \forall\, 0< r< +\infty. 
 \label{2.10}
 \end{equation}
 This shows that $ \varphi $ is nonnegative. Furthermore from \eqref{2.9} we get 
\[
\varphi ' =   \frac{1}{R}\intl_{\partial B_R} | \bu|^2  dS- \frac{\varphi }{R}=\frac{2}{R}\intl_{\partial B_R} | \bu _N|^2  dS\quad  \text{ a.\,e.  in $ (0, +\infty)$},  
\] 
and therefore $ \varphi $ is nondecreasing. In particular,
\begin{equation}
\exists\, \lim_{R \to 0^+} \varphi (R) \ge 0. 
\label{2.11}
\end{equation} 
\hspace{0.5cm}

{\bf Proof of Theorem\,\ref{thm1.3}} Observing \eqref{2.8} with $ \alpha =1$ along with $ \varphi (r) \ge 0$,  we infer 
\begin{equation}
2\intl_{B_R  \setminus B_r} \frac{| \bu _N|^2}{|x|} dx =\varphi (R) - \varphi (r) \le \varphi (R)
\quad \forall\, 0<r<R.  
 \label{2.14}
\end{equation}
 This shows that $ \frac{| \bu_N |^2}{| x|}\in  L^1_{ \rm loc}(\R^{3})$ and from monotone convergence of the measure 
 we conclude \eqref{2.12}.  

\hspace{0.5cm}
If the condition \eqref{2.13} is true,  then from \eqref{2.12} and $ \varphi (R) \le 
\intl_{\partial B_{ R}} | \bu_T |^2dS$ we get 
\[
\intl_{B_{ R_k}} \frac{| \bu_N |^2}{|x|}  dx \le \intl_{\partial B_{ R_k}} | \bu_T |^2 dS
\rightarrow 0\quad  \text{as}\quad  k \rightarrow +\infty. 
\]  
Accordingly, $ \frac{|\bu_N|^2}{|x|}=0$ which implies  $ \bu_N \equiv 0 $, and therefore 
$  \bu = \bu _T$. On the other hand from \eqref{2.8} with $ \alpha =1$ it follows that $ \varphi = \varphi _0 \equiv \const \ge 0$. 
Consequently, in view of \eqref{2.13} 
\[
\varphi _0 = \intl_{\partial B_{ R_k}} | \bu_T |^2dS \rightarrow 0 \quad  \text{as}\quad  k \rightarrow +\infty.
\]
Whence, $ \bu = \bu _T \equiv 0$.  \hfill \Beweisende \\
\hspace{0.5cm}

{\bf Proof of Theorem\,\ref{thm1.6}} First, by the monotonicity of $ \varphi  $ we observe 
\begin{equation}
\intl_{B_2} | \bu |^2  dx \ge \intl_{B_2} (| \bu_T |^2- | \bu_N |^2)  dx 
\ge  \intl_{1}^{2} \varphi (r)  dr \ge \varphi (1).   
\label{2.20}
\end{equation} 

By means of \eqref{2.1} we see that 
\begin{equation}
\frac{1}{R} \intl_{B_R} | \bu |^2  dx  \le \varphi (R) \le  \varphi (1)\quad  \forall\,0< R \le 1.     
\label{2.21}
\end{equation} 
Then the assertion of the theorem follows by combining  \eqref{2.20} with (\ref{2.21}).  \hfill \Beweisende  
\vspace{0.5cm}

$$\mbox{\bf Acknowledgements}$$
Chae was partially supported by NRF grant 2016R1A2B3011647, while Wolf has been supported by the Brain Pool Project of the Korea Federation of Science and Technology Societies  (141S-1-3-0022).  

\begin{thebibliography}{99}
\bibitem{arn} V. I. Arnold and B. A. Khesin, {\it Topological Methods in Hydrodynamics}, (1998), Springer-Verlag, New York.
\bibitem{cha0}D. Chae, {\it Conditions on the pressure for vanishing velocity in the incompressible fluid flows in $\Bbb R^N$,} Comm. PDE, {\bf 37,} (2012), pp. 1445-1455.
\bibitem{cha} D. Chae and P. Constantin, {\it Remarks on a Liouville-type theorem for Beltrami flows,}   Int. Math. Res. Notices,  {\bf 2015},  no.20,  (2015), pp. 10012-10016.
 \bibitem{con}P.  Constantin and A.  Majda,  {\it The Beltrami spectrum for incompressible fluid flows, } Comm. Math. Phys., {\bf 115},  no. 3, 
(1988), pp. 435-456. 
\bibitem{del}C. De Lellis and L. Sz\'{e}kelyhidi Jr., {\it Dissipative continuous Euler flows,} Invent Math., {\bf 193},  (2013), pp. 377-407.
\bibitem{enc}A. Enciso and D. Peralta-Salas,  {\it Knots and links in steady solutions of the Euler equation,} 
Ann. Math.,  {\bf 175}(2), no. 1 (2012), pp. 345-367.
\bibitem{nad}N. Nadirashvili, {\it Liouville theorem for Beltrami flow,} Geom. Funct. Anal., {\bf 24}, (2014), pp. 916-921.
\bibitem{pel} R. Pelz, V. Yakhot and S. A.  Orszag,  L. Shtilman, and E.  Levich, E., {\it Velocity-vorticity patterns in
turbulent flow, } Phys. Rev. Lett., {\bf 54},   (1985), pp. 2505-2509.
\end{thebibliography}
\end{document}